\documentclass[11pt]{amsart}
\usepackage[T1]{fontenc}
\usepackage{lmodern}
\usepackage{amsmath,amsthm,amssymb,amsfonts}
\usepackage{enumerate}
\usepackage{epsfig}
\usepackage[dvipsnames,svgnames,table]{xcolor}
\usepackage{pgf,tikz}
\usetikzlibrary{calc}
\usepackage{fullpage}
\usepackage[colorlinks=true,linkcolor=Blue,urlcolor=RoyalBlue,citecolor=red]{hyperref}
\usepackage{booktabs}
\usepackage[nameinlink]{cleveref}
\usepackage{enumitem}
\usepackage{mathtools}
\usepackage{faktor}
\usetikzlibrary{calc,shapes,positioning,patterns,arrows,decorations.pathreplacing}
\tikzset{Bullet/.style={fill=black,draw,color=#1,circle,minimum size=3pt,scale=0.75}}
\usepackage{xcolor}
\usepackage{mathdots}

\newtheorem{theorem}{Theorem}[section]
\newtheorem{proposition}[theorem]{Proposition}
\newtheorem{lemma}[theorem]{Lemma}
\newtheorem{corollary}[theorem]{Corollary}

\theoremstyle{definition}

\theoremstyle{remark}

\let\emph\relax
\DeclareTextFontCommand{\emph}{\bfseries\em}

\newcommand{\suchthat}{\;\ifnum\currentgrouptype=16 \middle\fi|\;}

\DeclareMathOperator{\sign}{sign}
\DeclareMathOperator{\Span}{Span}

\DeclareMathOperator{\rank}{rank}

\DeclareMathOperator{\conullity}{conullity}
\DeclareMathOperator{\coker}{Coker}
\DeclareMathOperator{\link}{Link}
\DeclareMathOperator{\nullity}{nullity}

\tikzstyle{vertex}=[fill=black,circle,inner sep=0pt, minimum size=4pt]
\tikzstyle{edge}=[line width=1.5pt,black]

\newcommand*{\defn}[1]{{\color{purple}\textit{#1}}}

\title{Counting frameworks of bipyramids}
\author{Jack Southgate}
\address{Jack Southgate, University of St Andrews}
\email{}
\urladdr{josouthgate.github.io}

\date{January 2026}

\begin{document}

\begin{abstract}
    We give a linear upper bound on the number of distinct volume-equivalent frameworks of bipyramids, up to rigid motions.
    As a corollary, we show that global volume rigidity is not a generic property of simplicial complexes.
\end{abstract}

\maketitle

\section{Introduction}\label{sec:intro}

Let $\Sigma$ be a $d$-dimensional abstract simplicial complex, assume throughout that $\Sigma$ is pure.
We may describe a realisation of $\Sigma$ in $\mathbb{R}^d$ by listing out the vectors of the positions of the vertices of $\Sigma$.
We will denote the individual vertices of $\Sigma$ as the numbers $1,\dots,n$ and their positions in $\mathbb{R}^d$ as $\mathbf{p}(1),\dots,\mathbf{p}(n)$.
The vector $\mathbf{p}=(\mathbf{p}(1),\dots,\mathbf{p}(n))\in(\mathbb{R}^d)^n$ is known as a \defn{configuration} and when it is paired with the additional combinatorial information of the rest of the simplicial complex, we obtain a \defn{framework}, $(\Sigma,\mathbf{p})$.

In \cref{sec:volrig}, we will introduce ways to measure the $d$-dimensional volumes of the $d$-simplices in frameworks in $\mathbb{R}^d$, and define a theory of volume rigidity.
Volume rigidity theory, although not as prolific as bar-joint rigidity, has been studied in several forms:
Tay et al. wrote about a slightly different form of volume rigidity to what we will consider here, motivated by algebraic combinatorics and $g$-theory \cite{tay1995skeletal1}, \cite{tay1995skeletal2}.
Streinu was the first to write about this specific problem from a rigidity theoretic point of view, with Theran \cite{streinu2007algorithms} and later Borcea \cite{borcea2012realizations}, \cite{borcea2014volume}.
Since then, Bulavka et al. used tools from algebraic combinatorics to study volume rigidity \cite{bulavka2022volume}, in particular, they prove a stronger version of \Cref{lem:vsrig} independently and by different means to this paper.
Finally, in the time between the upload of the initial preprint of this paper and now, Cruickshank, et al. \cite{cruickshank2025volume} and Lew et al. \cite{lew2025k} both uploaded preprints of work on the rigidity of frameworks of simplicial complexes in $\mathbb{R}^d$ where the volumes of lower-dimensional simplices are preserved.

Once we define volume rigidity theory, we are able to obtain an equivalence relation on $d$-dimensional frameworks of $(\Sigma,\mathbf{p})$ in $\mathbb{R}^d$, where two frameworks are equivalent if they are images of each other under a $d$-dimensional special affine transformation of $\mathbb{R}^d$ (ie. an affine transformation of $\mathbb{R}^d$ that preserves $d$-dimensional volumes).
The equivalence classes arising from this relation, known as \textit{congruence classes}, encode the distinct embeddings of $\Sigma$ up to rigid motions.

Borcea and Streinu began searching for bounds on the number of congruence classes using the degree of the \textit{measurement variety}, in bar-joint rigidity theory in \cite{borcea2002number} and applied similar methods to obtain the following bound for volume rigidity:

\begin{theorem}\label{thm:BSbound}\cite{borcea2012realizations}
    Let $\Sigma$ be a minimally volume rigid $d$-dimensional simplicial complex on $n\geq d+1$ vertices.
    A generic framework of $\Sigma$ admits at most
    \begin{equation}\label{eq:BSboundd}
        (d(n-d-1))!\prod\limits_{i=0}^{d-1}\frac{i!}{(n-d-1+i)!}
    \end{equation}
    congruence classes.
\end{theorem}

From \cref{sec:trisurf} onwards, we will consider the 2-dimensional case, where the bound \cref{eq:BSboundd} becomes
\begin{equation}\label{eq:BSbound2}
    \frac{1}{n-2}\binom{2n-6}{n-3},
\end{equation}
which is exponential in $n$, via Stirling's approximation.
Although its first few values do not diverge that quickly \cref{tab:BSbound2}, numerical experiments show that it is already an overestimate in these small cases when $\Sigma$ is a triangulation of $\mathbb{S}^2$ and even more so when $\Sigma$ is a bipyramid.
This motivates the main theorem of this paper.

\begin{theorem}\label{thm:main}
    Let $\Sigma$ be a bipyramid on $n\geq 5$ vertices.
    A generic framework of $\Sigma$ admits at most $n-4$ congruence classes.
\end{theorem}
\begin{centering}
\begin{table}[]
    \centering
    \begin{tabular}{|c||c|c|c|c|c|c|c|c|}
        \hline
        $n$ & 3 & 4 & 5 & 6 & 7 & 8 & 9 & 10 \\
        \hline
        $\frac{1}{n-2}\binom{2n-6}{n-3}$ & 1 & 1 & 2 & 5 & 14 & 42 & 132 & 429 \\
        \hline
    \end{tabular}
    \caption{First 8 values of \cref{eq:BSbound2}.}
    \label{tab:BSbound2}
\end{table}
\end{centering}

Global rigidity is a strengthening of rigidity theory which is well defined both in volume rigidity theory and bar-joint rigidity theory.
A framework is globally rigid if it admits only one congruence class.
A landmark result in bar-joint rigidity theory is the following theorem of Connelly and Gortler et al. stating that global bar-joint rigidity is a \textit{generic property} property of graphs.

\begin{theorem}\label{thm:CGHT}\cite{connelly2005generic},\cite{gortler2010characterizing}
    Let $G$ be a graph.
    Either all generic frameworks of $G$ in $\mathbb{R}^d$ are globally bar-joint rigid in $\mathbb{R}^d$ or none are.
\end{theorem}

By applying the methods used in the proof of \cref{thm:main}, we are able to show by counterexample that \cref{thm:CGHT} does not have a volume rigidity analogue.

\subsection{Acknowledgements}

This paper was written under the academic supervision of Louis Theran at the University of St Andrews.

The author would like to thank Louis Theran for advice and direction in developing the results of this paper and Alex Rutar for useful discussions.

\subsection{Notes on changes}

The first versions of this paper, available on the arXiv, were written in early 2023, while I was a PhD student.
They suffered from poor formatting and writing, as well as errors in some key proofs.
This version is a substantial rewrite of those initial versions.
No new original results have been added and some unnecessary results have been removed (they can now be found in my PhD thesis \cite{southgate2024volume} or in the work of Borcea and Streinu \cite{borcea2002number} or Bulavka et al. \cite{bulavka2022volume}).

This version of the paper also considers pure $d$-dimensional simplicial complexes instead of $(d+1)$-uniform simple hypergraphs.
This change was made to maintain consistency with much of the modern literature on volume rigidity.
The combinatorial data necessary for the results of this paper has been translated in a one-to-one way between these two objects.

\section{Volume rigidity}\label{sec:volrig}

In this section, we outline some basic definitions and results in volume rigidity theory.

As well as a vector, it is useful to think of configurations as corresponding to a matrix.
The \defn{configuration matrix} associated to the configuration $\mathbf{p}$ is the $(d+1)\times n$ matrix 
\begin{equation*}
    C(\mathbf{p})=\begin{bmatrix} 1 & \dots & 1 \\ \mathbf{p}(1) & \dots & \mathbf{p}(n) \end{bmatrix}.
\end{equation*}
For any $k$-tuple, written as the string $\sigma=i_1\dots i_k$, where $i_1<\dots<i_k$, $C(\sigma,\mathbf{p})$ will denote the submatrix of $C(\mathbf{p})$ whose columns are the columns $i_1,\dots,i_k$, in that order, of $C(\mathbf{p})$.

The \defn{complete measurement map} measures the volume of every $(d+1)$-tuple of vertices in a configuration and is defined as follows:
\begin{equation*}
    \alpha_n^d:(\mathbb{R}^d)^n\rightarrow\mathbb{R}^{\binom{n}{d+1}}; \mathbf{p}\mapsto\left(\det(C(\sigma,\mathbf{p})):\sigma\in\binom{[n]}{d+1}\right).
\end{equation*}
We will index the codomain, $\mathbb{R}^{\binom{n}{d+1}}$, of $\alpha_n^d$ by the $(d+1)$-tuples of $1,\dots,n$, ordered linearly lexicographically, denoting it $\mathbb{R}^{\binom{[n]}{d+1}}$.
By doing so, we are able to orthogonally project onto the coordinates indexed by the maximal simplices of a $d$-dimensional simplicial complex, denoting this projection $\pi_{\Sigma}$.
We are therefore able to define the \defn{measurement map} of a simplicial complex (on at most $n$ vertices) as $\alpha_{\Sigma}=\pi_{\Sigma}\circ\alpha_n^d$.

We will use $K_n^d$ to denote the complete $d$-dimensional simplicial complex on $n$ vertices, then $\alpha_{K_n^d}=\alpha_n^d$.

In this paper we will only be considering $d$-dimensional frameworks in $\mathbb{R}^d$, ie. frameworks $(\Sigma,\mathbf{p})$ for which the affine span of the configuration $\mathbf{p}$ is $d$-dimensional.
This condition is automatically met when a framework is in \defn{general position}, i.e. the span of each $k$-simplex is $k$-dimensional, or is \defn{generic}, i.e. the entries of the vector $\mathbf{p}$ are algebraically independent over $\mathbb{Q}$.

As a result, we can compare frameworks by their maximal simplices' volumes: two frameworks $(\Sigma,\mathbf{p})$ and $(\Sigma,\mathbf{q})$ are \defn{equivalent} if $\alpha_{\Sigma}(\mathbf{p})=\alpha_{\Sigma}(\mathbf{q})$.
The complete measurement map is the map $\alpha_n^d=\alpha_{K_n^d}$, where $K_n^d$ is the completed pure $d$-dimensional simplicial complex on $n$ vertices, it measures the $d$-dimensional volumes of all $(d+1)$-tuples in a simplicial complex.
Two frameworks $(\Sigma,\mathbf{p})$ and $(\Sigma,\mathbf{q})$ are \defn{congruent} if $\alpha_n^d(\mathbf{p})=\alpha_n^d(\mathbf{q})$.

We are now in a position to define ($d$-dimensional volume) rigidity (in $\mathbb{R}^d$).
As with the bar-joint case, there are many equivalent definitions of rigidity, this first definition is analogous to the most commonly stated definition from bar-joint rigidity theory.

We say that a framework $(\Sigma,\mathbf{p})$ is \defn{rigid} if there exists an open subset $U$ of $\mathbf{p}$ in $(\mathbb{R}^d)^n$ so that, for every $q\in U$, if $(\Sigma,\mathbf{p})$ and $(\Sigma,\mathbf{q})$ are equivalent, then they are congruent.
If $(\Sigma,\mathbf{p})$ fails to be rigid, then it is \defn{flexible}.

A flexible framework therefore admits a continuous family of deformations that remain equivalent to the original framework, whilst a rigid framework does not.
This idea is formalised in the language of flexes.

A \defn{finite flex} of the framework $(\Sigma,\mathbf{p})$ in $\mathbb{R}^d$ is a continuous (in the Euclidean topology) map $\gamma:[0,1]\rightarrow(\mathbb{R}^d)^n$ such that:
\begin{itemize}
    \item[F1]\label{item:f1} $\gamma(0)=\mathbf{p}$;
    \item[F2]\label{item:f2} $(\Sigma,\gamma(t))$ is equivalent to $(\Sigma,\mathbf{p})$, for each $t\in [0,1]$.
\end{itemize}
A finite flex is \defn{trivial} if it meets the additional condition
\begin{itemize}
    \item[F3]\label{item:f3} $(\Sigma,\gamma(t))$ is congruent to $(\Sigma,\mathbf{p})$, for each $t\in [0,1]$.
\end{itemize}
With this language, we are able to provide a second definition of rigidity.

The framework $(\Sigma,\mathbf{p})$ is rigid if every finite flex that $(\Sigma,\mathbf{p})$ admits is trivial.

\begin{proposition}\label{prop:rigeq}
    The two definitions of rigidity given above are equivalent.
\end{proposition}

In order to prove \cref{prop:rigeq}, we will require the curve-selection Lemma, stated below, a standard result from real analysis.

\begin{lemma}\label{lem:CSL}\cite{milnor}
    Let $S\subseteq\mathbb{R}^D$ be a semi-algebraic set.
    Let $\mathbf{x}\in S$ and let $U$ be an open neighbourhood of $\mathbf{x}$ in $\mathbb{R}^D$.
    For each $\mathbf{y}\in S\cap U$, there exists an analytic semi-algebraic curve $\gamma:[0,1]\rightarrow S\cap U$ with $\gamma(0)=\mathbf{x}$ and $\gamma(1)=\mathbf{y}$.
\end{lemma}

\begin{proof}[Proof of \cref{prop:rigeq}]
    Let $(\Sigma,\mathbf{p})$ be a framework and assume that every finite flex of $(\Sigma,\mathbf{p})$ is trivial.
    Suppose, for the sake of contradiction, that every open neighbourhood of $\mathbf{p}$ in $(\mathbb{R}^d)^n$ contains a configuration yielding a framework to which $(\Sigma,\mathbf{p})$ is equivalent, but not congruent.
    Take any such framework, $(\Sigma,\mathbf{q})$, then, by \Cref{lem:CSL}, there exists a finite flex from $(\Sigma,\mathbf{p})$ to $(\Sigma,\mathbf{q})$.
    Since $(\Sigma,\mathbf{p})$ is not congruent to $(\Sigma,\mathbf{q})$, this finite flex is non-trivial, a contradiction.

    Let $(\Sigma,\mathbf{p})$ be a framework and assume that there exists an open neighbourhood $U$ of $\mathbf{p}$ in $(\mathbb{R}^d)^n$ such that if $q\in U$ and $(\Sigma,\mathbf{p})$ and $(\Sigma,\mathbf{q})$ are equivalent, then $(\Sigma,\mathbf{p})$ and $(\Sigma,\mathbf{q})$ are congruent.
    Suppose that $\gamma$ is a non-trivial finite flex of $(\Sigma,\mathbf{p})$, then by its continuity, there exists $\tau\in[0,1)$ so that $\gamma(\tau)\in U$, but $(\Sigma,\gamma(\tau))$ is equivalent but not congruent to $(\Sigma,\mathbf{p})$, a contradiction.
\end{proof}

The \defn{rigidity matrix} of a framework $(\Sigma,\mathbf{p})$, denoted $R(\Sigma,\mathbf{p})$, is the $(f(\Sigma)_d\times dn)$-matrix obtained by evaluating the Jacobian of $\alpha_{\Sigma}$ at $\mathbf{p}$.
An \defn{infinitesimal flex} of $(\Sigma,\mathbf{p})$ is the infinitesimal velocity of a finite flex of $(\Sigma,\mathbf{p})$.
An infinitesimal flex of $(\Sigma,\mathbf{p})$ is trivial if it is the infinitesimal velocity of a trivial finite flex and non-trivial otherwise.

\begin{proposition}\label{prop:rigker}\cite[p. 48]{southgate2024volume}
    The kernel vectors of $R(\Sigma,\mathbf{p})$ are precisely the infinitesimal flexes of $(\Sigma,\mathbf{p})$.
\end{proposition}

Let $\Sigma$ be a simplicial complex on $n$ vertices.
We say that the framework $(\Sigma,\mathbf{p})$ in $\mathbb{R}^d$ is \defn{infinitesimally rigid} if the only infinitesimal flexes of $(\Sigma,\mathbf{p})$ are trivial.
By \cref{prop:rigker}, $(\Sigma,\mathbf{p})$ being infinitesimally rigid is equivalent to its rigidity matrix having rank $dn-(d^2+d-1)$.
Note that infinitesimal rigidity is not equivalent to rigidity, as seen in \cref{fig:octinfflex}.
However, these notions of rigidity are equivalent in the case of generic frameworks.

\begin{figure}
        \centering
        \begin{tikzpicture}
            \filldraw[fill=yellow!80!black,line width=0.5pt,fill opacity=0.6] (0,0) -- (3,0) -- (0,3) -- (0,0);
            \filldraw[fill=yellow!80!black,line width=0.5pt,fill opacity=0.6] (0,0) -- (3,0) -- (1.4,0.5) -- (0,0);
            \filldraw[fill=yellow!80!black,line width=0.5pt,fill opacity=0.6] (0,0) -- (0,3) -- (0.5,1.4) -- (0,0);
            \filldraw[fill=yellow!80!black,line width=0.5pt,fill opacity=0.6] (3,0) -- (0,3) -- (1.2,1.2) -- (3,0);
            \filldraw[fill=yellow!80!black,line width=0.5pt,fill opacity=0.6] (3,0) -- (1.4,0.5) -- (1.2,1.2) -- (3,0);
            \filldraw[fill=yellow!80!black,line width=0.5pt,fill opacity=0.6] (0,3) -- (0.5,1.4) -- (1.2,1.2) -- (0,3);
            \filldraw[fill=yellow!80!black,line width=0.5pt,fill opacity=0.6] (0,0) -- (1.4,0.5) -- (0.5,1.4) -- (0,0);
            \filldraw[fill=yellow!80!black,line width=0.5pt,fill opacity=0.6] (1.4,0.5) -- (0.5,1.4) -- (1.2,1.2) -- (1.4,0.5);
            \coordinate[Bullet=black, label=below:1] (n1) at (0,0);
            \coordinate[Bullet=black, label=below:2] (n2) at (3,0);
            \coordinate[Bullet=black, label=left:3] (n3) at (0,3);
            \coordinate[Bullet=black, label=left:{4}] (n4) at (0.5,1.4);
            \coordinate[Bullet=black, label=below:{5}] (n5) at (1.4,0.5);
            \coordinate[Bullet=black, label=right:{6}] (n6) at (1.2,1.2);
            \draw[->,red] (1.4,0.5) -- (2,0.5);
            \draw[->,red] (0.5,1.4) -- (0.5,0.8);
            \draw[->,red] (1.2,1.2) -- (0.8,1.6);
            \draw[dashed] (-0.2,-0.2) -- (1.7,1.7);
        \end{tikzpicture}
        \caption{The octahedron is generically rigid (as we will show in \cref{prop:rigs2}), but frameworks of the octahedron that are mirror-symmetric in the line $\Span\{\mathbf{p}(0),\mathbf{p}(6)\}$ admit an infinitesimal flex.}
        \label{fig:octinfflex}
\end{figure}

\begin{proposition}\cite{borcea2012realizations}\cite[p. 50]{southgate2024volume}
When $\mathbf{p}$ is generic, $(\Sigma,\mathbf{p})$ is generic if and only if $\rank(R(\Sigma,\mathbf{p}))=dn-(d^2+d-1)$.

Rigidity is a generic property of $\Sigma$, ie. either all generic frameworks of $\Sigma$ are rigid or none are.
\end{proposition}

Since rigidity is a generic property of simplicial complexes, it makes sense to talk about rigidity in terms of simplicial complexes instead of just their frameworks.
A simplicial complex $\Sigma$ is \defn{rigid} if all of its generic frameworks are rigid, and \defn{flexible} otherwise.
A simplicial complex $\Sigma$ is \defn{minimally rigid} if it is generic, but any simplicial complex $\Sigma'$, defined by $\Sigma'^{(d)}=\Sigma^{(d)}\setminus\{\sigma\}$, is flexible.
A simplicial complex $\Sigma$ is \defn{redundantly rigid} if it is rigid and every simplicial complex $\Sigma'$ as above is rigid.

\section{Counting frameworks}

Both flexible and rigid frameworks in $\mathbb{R}^d$ will admit uncountably infinitely many equivalent frameworks.
Indeed, all frameworks admit their images under special affine transformations of $\mathbb{R}^d$ (ie. trivial finite flexes), while flexible ones will additionally admit their images under non-trivial finite flexes.

In both cases, however, we can ignore images under trivial finite flexes by considering the quotient space associated to each framework $(\Sigma,\mathbf{p})$:
\begin{equation*}
    \mathcal{C}(\Sigma,\mathbf{p}) = \faktor{\alpha_{\Sigma}^{-1}(\alpha_{\Sigma}(\mathbf{p}))}{(\alpha_n^d)^{-1}(\alpha_n^d(\mathbf{p}))}.
\end{equation*}
We call $\mathcal{C}(\Sigma,\mathbf{p})$ the \defn{configuration space} of $(\Sigma,\mathbf{p})$ and its elements, $[q]$, \defn{congruence classes}, as the quotient above defines an equivalence relation under congruence.

\begin{lemma}
    For each $(\Sigma,\mathbf{p})$, $\mathcal{C}(\Sigma,\mathbf{p})$ has finitely many connected components.
\end{lemma}

\begin{proof}
    We may think of $\mathcal{C}(\Sigma,\mathbf{p})$ as the quotient of the semi-algebraic variety $V = \alpha_{\Sigma}^{-1}(\alpha_{\Sigma}(\mathbf{p}))$ by the action of the special affine group $\mathrm{SA}(d,\mathbb{R})$.
    Since $V$ is semi-algebraic, it has finitely many connected components (see most introductory texts on real algebraic geometry, for example \cite{bochnak2013real}), and therefore so does its quotient.
\end{proof}

\begin{lemma}
    The configuration space $\mathcal{C}(\Sigma,\mathbf{p})$ is zero-dimensional if and only if $(\Sigma,\mathbf{p})$ is rigid.
\end{lemma}

\begin{proof}
    Suppose that $(\Sigma,\mathbf{p})$ is flexible.
    Then, within any open neighbourhood $U$ of $\mathbf{p}$ in $\mathcal{C}(\Sigma,\mathbf{p})$, equipped with the subspace topology, there exists $\mathbf{q}\in U\setminus\{\mathbf{p}\}$ and, by \Cref{lem:CSL}, a continuous path $\gamma:[0,1]\rightarrow U$ such that $\gamma(0)=\mathbf{p}$ and $\gamma(1)=\mathbf{q}$.
    Hence $\mathcal{C}(\Sigma,\mathbf{p})$ is at least one-dimensional.

    Suppose that $\mathcal{C}(\Sigma,\mathbf{p})$ is zero-dimensional.
    Let $C$ be a connected component of $\mathcal{C}(\Sigma,\mathbf{p})$.
    Let $\mathbf{q}_1\in C$ be generic and let $U$ be an open neighbourhood of $\mathbf{q}_1$, in $\mathcal{C}(\Sigma,\mathbf{p})$, equipped with the subspace topology.
    Let $\mathbf{q}_2\in U$.
    If there were a continuous, non-constant, path from $\mathbf{q}_1$ to $\mathbf{q}_2$ in $U$, this would constitute a one-dimensional subspace of $\mathcal{C}(\Sigma,\mathbf{p})$, hence $\mathbf{q}_2=\mathbf{q}_1$ and $(\Sigma,\mathbf{q}_1)$ is rigid.
    Since rigidity is a generic property, $(\Sigma,\mathbf{p})$ is rigid too.
\end{proof}

The number of configuration classes of different frameworks of the same simplicial complex in $\mathbb{R}^d$ is not necessarily an invariant of $\Sigma$, even when restricting to generic frameworks, as we will see in \cref{cor:notggr}.
Let $c(\Sigma,\mathbf{p})$ denote the number of configuration classes of $(\Sigma,\mathbf{p})$ and let $c(\Sigma) = \max\limits_{\substack{\mathbf{p}\in(\mathbb{R}^d)^n \\ \mathbf{p}\text{ generic}}}\{c(\Sigma,\mathbf{p})\}$.

We now restate \cref{thm:BSbound} in using our new notation.

\begin{theorem}\cite{borcea2012realizations}
    Let $\Sigma$ be a minimally rigid $d$-dimensional simplicial complex on $n$ variables, then $c(\Sigma) \leq (d(n-d-1))!\prod\limits_{i=0}^{d-1}\frac{i!}{(n-d-1+1)!}$.
\end{theorem}

Combinatorial operations on simplicial complexes can affect the number of congruence classes they admit in predictable ways.

Let $\Sigma_1$ and $\Sigma_2$ be two $d$-dimensional simplicial complexes and let $\sigma^1\in\Sigma_1^{(d)}$ and $\sigma^2\in\Sigma_2^{(d)}$.
We may \defn{glue} the two complexes to obtain $\Sigma_1*\Sigma_2$ by identifying $\sigma_i^1$ with $\sigma_{\rho(i)}$, for each $i\in[d+1]$, with $\rho\in S_{d+1}$ a permutation.

Let $(\Sigma_1,\mathbf{p}_1)$ and $(\Sigma_2,\mathbf{p}_2)$ be two frameworks in $\mathbb{R}^d$ and suppose that $\mathbf{p}_1(\sigma_i^1)=\mathbf{p}_2(\sigma_{\rho(i)}^2)$, for each $i\in [d+1]$.
We may glue the two frameworks to obtain $(\Sigma_1*\Sigma_2,\mathbf{p}_1*\mathbf{p}_2)$ by setting
\begin{equation*}
    (\mathbf{p}_1*\mathbf{p}_2)(i) = \begin{cases}
        \mathbf{p}_1(i),&\text{ if }i\in\Sigma_1^{(0)},\\
        \mathbf{p}_2(i),&\text{ if }i\in\Sigma_2^{(0)}.
    \end{cases}
\end{equation*}

\begin{proposition}\label{prop:gluing}
    Let $(\Sigma_1,\mathbf{p}_1)$ and $(\Sigma_2,\mathbf{p}_2)$ be two frameworks as above.
    If both $(\Sigma_1,\mathbf{p}_1)$ and $(\Sigma_2,\mathbf{p}_2)$ are rigid, then so is $(\Sigma_1*\Sigma_2,\mathbf{p}_1*\mathbf{p}_2)$ and moreover $c(\Sigma_1*\Sigma_2,\mathbf{p}_1*\mathbf{p}_2) = c(\Sigma_1,\mathbf{p}_1)c(\Sigma_2,\mathbf{p}_2)$.
\end{proposition}

In proving \cref{prop:gluing}, we note the following equivalent definitions of equivalence and congruence of frameworks.
Two frameworks $(\Sigma,\mathbf{p})$ and $(\Sigma,\mathbf{q})$ are equivalent if, for each $\sigma\in\Sigma^{(d)}$, there exists a special affine transformation $f_{\sigma}$ of $\mathbb{R}^d$ such that $f_{\sigma}(\mathbf{p}(i))=\mathbf{q}(i)$, for each $i\in\sigma^{(0)}$.
They are congruent if, moreover, there exists a single special affine transformation $f$ of $\mathbb{R}^d$ such that $f(\mathbf{p}(i))=\mathbf{q}(i)$, for each $i\in\Sigma^{(0)}$.

The equivalence of these definitions with those given in \cref{sec:volrig} can be seen by the fact that an affine transformation of $\mathbb{R}^d$ is uniquely defined by its action on any set of $d+1$ affinely independent points.

\begin{proof}
    Write $\Sigma=\Sigma_1*\Sigma_2$, $p=\mathbf{p}_1*\mathbf{p}_2$, $n_1 = f(\Sigma_1)_0$, $n_2=f(\Sigma_2)_0$, $n=n_1+n_2-d-1=f(\Sigma)_0$.
    
    Let $U_1$ (resp. $U_2$) be the neighbourhood of $\mathbf{p}_1$ (resp. $\mathbf{p}_2$) in $(\mathbb{R}^d)^{n_1}$ (resp. $(\mathbb{R}^d)^{n_2}$) on which equivalence implies congruence.
    Let $U\subseteq\mathbf{q}(\mathbb{R}^d)^n$ be such that $\pi_{\Sigma_1^{(0)}}(U)=U_1$ and $\pi_{\Sigma_2^{(0)}}(U)=U_2$.
    Let $q\in U$ and suppose that $(\Sigma,\mathbf{p})$ and $(\Sigma,\mathbf{q})$ are equivalent.
    Write $\mathbf{q}_1=\pi_{\Sigma_1^{(0)}}(\mathbf{q})$ and $\mathbf{q}_2=\pi_{\Sigma_1^{(0)}}(\mathbf{q})$, then $\mathbf{q}_1\in U_1$ and $\mathbf{q}_2\in U_2$, so, by the rigidity of $\Sigma_1$ (resp. $\Sigma_2$), $(\Sigma_1,\mathbf{p}_1)$ (resp. $(\Sigma_2,\mathbf{p}_2)$) is congruent to $(\Sigma_1,\mathbf{q}_1)$ (resp. $(\Sigma_2,\mathbf{q}_2)$).
    Let $f_1$ (resp. $f_2$) be the affine transformations sending each vertex of $(\Sigma_1,\mathbf{p}_1)$ (resp. $(\Sigma_2,\mathbf{p}_2)$) to its position in $(\Sigma_1,\mathbf{q}_1)$ (resp. $(\Sigma_2,\mathbf{q}_2)$).
    Then $f_1$ and $f_2$ agree on their action on the $d+1$ vertices of $\sigma$, hence they must be equal.

    Next, suppose that $(\Sigma_1,\mathbf{p}_1)$ is equivalent, but not congruent to $(\Sigma_1,\mathbf{q}_1)$ with, after applying a special affine transformation to the whole framework if necessary, $\mathbf{q}_1(\sigma_i^1)=\mathbf{p}_1(\sigma_i^1)$, for all $i\in [d+1]$.
    Then $(\Sigma,\mathbf{p})$ is equivalent, but not congruent to $(\Sigma,\mathbf{q}_1*\mathbf{p}_2)$.
    Applying the same argument to equivalent frameworks of $(\Sigma_2,\mathbf{p}_2)$, we see that $c(\Sigma,\mathbf{p})\geq c(\Sigma_1,\mathbf{p}_1)c(\Sigma_2,\mathbf{p}_2)$.

    Now, suppose that $(\Sigma,\mathbf{p})$ is equivalent, but not congruent to $(\Sigma,\mathbf{q})$.
    Then, by applying the same projections as in the second paragraph, we see that $(\Sigma_1,\mathbf{p}_1)$ (resp. $(\Sigma_2,\mathbf{p}_2)$) is equivalent to $(\Sigma_1,\mathbf{q}_1)$ (resp. $(\Sigma_2,\mathbf{q}_2)$).
    Therefore $c(\Sigma,\mathbf{p})\leq c(\Sigma_1,\mathbf{p}_1)c(\Sigma_2,\mathbf{p}_2)$, concluding the proof.
\end{proof}

A framework is \defn{globally rigid} in $\mathbb{R}^d$ if it only admits a single congruence class.

\begin{proposition}\label{prop:compgr}
    Let $(K_n^d,\mathbf{p})$ be a generic framework in $\mathbb{R}^d$, with $n\geq d+1$, then $(K_n^d,\mathbf{p})$ is globally rigid.
\end{proposition}

\begin{proof}
    Firstly, $(K_{d+1}^d,\mathbf{p})$ is globally rigid.
    Indeed, since it consists of a single $d+1$ simplex any two equivalent frameworks are immediately also congruent.

    Now, suppose that $(K_n^d,\mathbf{p})$ is globally rigid, for some $n\geq d+1$.
    Let $(K_{n+1}^d,\mathbf{p}^*)$ be a generic framework with $\mathbf{p}^*(i)=\mathbf{p}(i)$, for all $i\in[n]$.
    Suppose that $(K_{n+1}^d,\mathbf{p}^*)$ is equivalent to $(K_{n+1}^d,\mathbf{q})$, then $(K_n^d,\mathbf{p})$ is equivalent, and therefore congruent to $(K_n^d,\pi_{[n]}(\mathbf{q}))$.
    We can therefore apply a special affine transformation to $(K_{n+1}^d,\mathbf{q})$ to obtain $(K_{n+1}^d,\mathbf{q}')$, with $\mathbf{q}'(i)=\mathbf{p}(i)$, for all $i\in [n]$.
    We claim that the position of $\mathbf{q}'(n+1)$ is uniquely defined as $\mathbf{p}(n+1)$.

    The position of $\mathbf{q}'(n+1)$ is defined as the intersection of $\binom{n}{d}>d$ hyperplanes in $\mathbb{R}^d$ (ie. either a single point, if it exists, or nothing), each parallel to $\Span\left\{\mathbf{q}'(i):i\in\binom{[n]}{d}\right\}=\Span\left\{\mathbf{p}(i):i\in\binom{[n]}{d}\right\}$ and at a distance of $\frac{1}{d!}\det(C(\tau (n+1),\mathbf{q}'))=\frac{1}{d!}\det(C(\tau(n+1),\mathbf{p}))$.
    Moreover, such an intersection must exist, as it is realised by $\mathbf{p}(n+1)$.
    Therefore $\mathbf{q}'(n+1)=\mathbf{p}(n+1)$, so $(K_{n+1}^d,\mathbf{q}')=(K_{n+1}^d,\mathbf{p})$, and hence $(K_{n+1}^d,\mathbf{p})$ and $(K_{n+1}^d,\mathbf{q})$ are congruent.
\end{proof}

The proof of \cref{prop:compgr} hints at a technique that we will use in counting congruence classes going forward: \textit{pinning} frameworks.
We \defn{pin} a framework $(\Sigma,\mathbf{p})$ in $\mathbb{R}^d$ by fixing the position of one of its $d$-simplices, and that of all frameworks equivalent to it.
In doing so, we mod $\alpha_{\Sigma}^{-1}(\alpha_{\Sigma}(\mathbf{p}))$ out by special affine transformations of $\mathbb{R}^d$, by composing them with the special affine transformation taking the vertices of the pinned simplex back to their positions in $(\Sigma,\mathbf{p})$.
We denote the pinned framework of $(\Sigma,\mathbf{p})$ by $(\Sigma,\overline{\mathbf{p}})$.

As in the proof of \cref{prop:compgr}, two pinned frameworks of $\Sigma$ are congruent if and only if they are equal \cite[p. 56]{southgate2024volume}.

\section{Triangulations of surfaces}\label{sec:trisurf}

Let $\Sigma$ be a pure 2-dimensional simplicial complex and let $u\in\Sigma^{(0)}$.
A \defn{vertex split} of $\Sigma$ at $u$ deletes $m$ 2-simplices containing $u$: $u\tau_1,\dots,u\tau_m$ and adds a new vertex $u^*$ as well as $m+2$ new 2-simplices $uu^*\tau_1^1, uu^*\tau_m^2,u^*\tau_1,\dots,u^*\tau_m$, where each $\tau_i=\tau_i^1\tau_i^2$, with $\tau_i^1<\tau_i^2$.
\Cref{fig:vertexsplit} demonstrates the vertex splitting process.

\begin{figure}
    \centering
        \begin{tikzpicture}
            \filldraw[fill=yellow!80!black,line width=0.5pt,fill opacity=0.3] (0,0) -- (0,2) -- (1.414,1.414) -- (0,0);
            \filldraw[fill=yellow!80!black,line width=0.5pt,fill opacity=0.3] (0,0) -- (1.414,1.414) -- (2,0) -- (0,0);
            \filldraw[fill=yellow!80!black,line width=0.5pt,fill opacity=0.6] (0,0) -- (1.414,1.414) -- (2,0) -- (0,0);
            \filldraw[fill=yellow!80!black,line width=0.5pt,fill opacity=0.3] (0,0) -- (0,-2) -- (-1.414,-1.414) -- (0,0);
            \filldraw[fill=yellow!80!black,line width=0.5pt,fill opacity=0.3] (0,0) -- (-2,0) -- (-1.414,-1.414) -- (0,0);
            \filldraw[fill=yellow!80!black,line width=0.5pt,fill opacity=0.6] (0,0) -- (-2,0) -- (-1.414,-1.414) -- (0,0);
            \coordinate[Bullet=black, label=left:{$u$}] (n1) at (0,0);
            \coordinate[Bullet=black, label=above:{$\tau_1^1$}] (n2) at (0,2);
            \coordinate[Bullet=black, label=right:$\tau_1^2$] (n3) at (1.414,1.414);
            \coordinate[Bullet=black, label=right:$\tau_{m_1}^2$] (n4) at (2,0);
            \coordinate[Bullet=black, label=below:$\tau_{m_k}^1$] (n6) at (0,-2);
            \coordinate[Bullet=black, label=left:$\tau_{m_k}^2$] (n7) at (-1.414,-1.414);
            \coordinate[Bullet=black, label=left:$\tau_m^2$] (n8) at (-2,0);
            \coordinate[label=right:{$\iddots$}] (n11) at (1,-1);
            \coordinate[label=right:{$\iddots$}] (n12) at (-1,1);
        \end{tikzpicture}
        \begin{tikzpicture}
            \filldraw[fill=yellow!80!black,line width=0.5pt,fill opacity=0.3] (0.5,-1) -- (0,2) -- (1.414,1.414) -- (0.5,-1);
            \filldraw[fill=yellow!80!black,line width=0.5pt,fill opacity=0.3] (0.5,-1) -- (-1.414,1.414) -- (-0.5,1) -- (0.5,-1);
            \filldraw[fill=yellow!80!black,line width=0.5pt,fill opacity=0.3] (0.5,-1) -- (1.414,1.414) -- (2,0) -- (0.5,-1);
            \filldraw[fill=yellow!80!black,line width=0.5pt,fill opacity=0.6] (0.5,-1) -- (1.414,1.414) -- (2,0) -- (0.5,-1);
            \filldraw[fill=yellow!80!black,line width=0.5pt,fill opacity=0.3] (0.5,-1) -- (0,-2) -- (-1.414,-1.414) -- (0.5,-1);
            \filldraw[fill=yellow!80!black,line width=0.5pt,fill opacity=0.3] (0.5,-1) -- (-2,0) -- (-1.414,-1.414) -- (0.5,-1);
            \filldraw[fill=yellow!80!black,line width=0.5pt,fill opacity=0.6] (0.5,-1) -- (-2,0) -- (-1.414,-1.414) -- (0.5,-1);
            \filldraw[fill=yellow!80!black,line width=0.5pt,fill opacity=0.3] (0.5,-1) -- (-2,0) -- (-1.414,1.414) -- (0.5,-1);
            \filldraw[fill=yellow!80!black,line width=0.5pt,fill opacity=0.3] (0.5,-1) -- (-0.5,1) -- (0,2) -- (0.5,-1);
            \coordinate[Bullet=black, label=above:{$u$}] (n1) at (-0.5,1);
            \coordinate[Bullet=black, label=right:{$u^*$}] (n10) at (0.5,-1);
            \coordinate[Bullet=black, label=above:{$\tau_1^1$}] (n2) at (0,2);
            \coordinate[Bullet=black, label=right:$\tau_1^2$] (n3) at (1.414,1.414);
            \coordinate[Bullet=black, label=right:$\tau_{m_1}^2$] (n4) at (2,0);
            \coordinate[Bullet=black, label=left:$\tau_{m_k}^1$] (n6) at (0,-2);
            \coordinate[Bullet=black, label=left:$\tau_{m_k}^2$] (n7) at (-1.414,-1.414);
            \coordinate[Bullet=black, label=left:$\tau_m^1$] (n8) at (-2,0);
            \coordinate[Bullet=black, label=left:$\tau_m^2$] (n9) at (-1.414,1.414);
            \coordinate[label=right:{$\iddots$}] (n11) at (1,-1);
        \end{tikzpicture}
    \caption{The vertex splitting process, with darker sections corresponding to multiple 2-simplices}
    \label{fig:vertexsplit}
\end{figure}

\begin{figure}
        \centering
        \begin{tikzpicture}
            \fill[yellow!80!, opacity = 0.3] (0,0) circle[radius=22.5mm];
            \filldraw[fill=yellow!80!black,line width=0.5pt,fill opacity=0.6] (0,0) -- (1.31,0.75) -- (1.31,-0.75) -- (0,0);
            \filldraw[fill=yellow!80!black,line width=0.5pt,fill opacity=0.6] (0,0) -- (-1.31,-0.75) -- (-1.31,0.75) -- (0,0);
            \coordinate[Bullet=black, label=left:{$u$}] (n1) at (0,0);
            \coordinate[Bullet=black, label=right:{$v_1$}] (n2) at (0,1.5);
            \coordinate[Bullet=black, label=right:$v_2$] (n3) at (1.31,0.75);
            \coordinate[Bullet=black, label=right:$v_k$] (n4) at (1.31,-0.75);
            \coordinate[Bullet=black, label=left:$w_1$] (n5) at (0,-1.5);
            \coordinate[Bullet=black, label=left:$w_2$] (n6) at (-1.31,-0.75);
            \coordinate[Bullet=black, label=left:$w_{\ell}$] (n7) at (-1.31,0.75);
            \coordinate[label=left:$\vdots$] (n9) at (-1.31,0);
            \coordinate[label=right:$\vdots$] (n10) at (1.31,0);
            \draw[-] (n1) -- (n2);
            \draw[-] (n1) -- (n3);
            \draw[-] (n1) -- (n4);
            \draw[-] (n1) -- (n5);
            \draw[-] (n1) -- (n6);
            \draw[-] (n1) -- (n7);
            \draw[-] (n2) -- (n3);
            \draw[-] (n2) -- (n7);
            \draw[-] (n3) -- (n4);
            \draw[-] (n4) -- (n5);
            \draw[-] (n5) -- (n6);
            \draw[-] (n6) -- (n7);
        \end{tikzpicture}
        \begin{tikzpicture}
            \fill[yellow!80!, opacity = 0.3] (0,0) circle[radius=22.5mm];
            \filldraw[fill=white] (0,0) -- (0,1.5) -- (1.31,0.75) -- (1.31,-0.75) -- (0,0);
            \filldraw[fill=yellow!80!black,line width=0.5pt,fill opacity=0.6] (0,0) -- (-1.31,-0.75) -- (-1.31,0.75) -- (0,0);
            \coordinate[Bullet=black, label=left:{$u$}] (n1) at (0,0);
            \coordinate[Bullet=black, label=right:{$v_1$}] (n2) at (0,1.5);
            \coordinate[Bullet=black, label=right:$v_2$] (n3) at (1.31,0.75);
            \coordinate[Bullet=black, label=right:$v_k$] (n4) at (1.31,-0.75);
            \coordinate[Bullet=black, label=left:$w_1$] (n5) at (0,-1.5);
            \coordinate[Bullet=black, label=left:$w_2$] (n6) at (-1.31,-0.75);
            \coordinate[Bullet=black, label=left:$w_{\ell}$] (n7) at (-1.31,0.75);
            \coordinate[label=left:$\vdots$] (n9) at (-1.31,0);
            \coordinate[label=right:$\vdots$] (n10) at (1.31,0);
            \draw[-] (n1) -- (n5);
            \draw[-] (n1) -- (n6);
            \draw[-] (n1) -- (n7);
            \draw[-] (n2) -- (n3);
            \draw[-] (n2) -- (n7);
            \draw[-] (n3) -- (n4);
            \draw[-] (n4) -- (n5);
            \draw[-] (n5) -- (n6);
            \draw[-] (n6) -- (n7);
        \end{tikzpicture}
        \begin{tikzpicture}
            \fill[yellow!80!, opacity = 0.3] (0,0) circle[radius=22.5mm];
            \filldraw[fill=yellow!80!black,line width=0.5pt,fill opacity=0.6] (0.65,0.38) -- (1.31,-0.75) -- (1.31,0.75) -- (0.65,0.38);
            \filldraw[fill=yellow!80!black,line width=0.5pt,fill opacity=0.6] (0,0) -- (-1.31,-0.75) -- (-1.31,0.75) -- (0,0);
            \coordinate[Bullet=black, label=left:{$u$}] (n1) at (0,0);
            \coordinate[Bullet=black, label=right:{$v_1$}] (n2) at (0,1.5);
            \coordinate[Bullet=black, label=right:$v_2$] (n3) at (1.31,0.75);
            \coordinate[Bullet=black, label=right:$v_k$] (n4) at (1.31,-0.75);
            \coordinate[Bullet=black, label=left:$w_1$] (n5) at (0,-1.5);
            \coordinate[Bullet=black, label=left:$w_2$] (n6) at (-1.31,-0.75);
            \coordinate[Bullet=black, label=left:$w_{\ell}$] (n7) at (-1.31,0.75);
            \coordinate[Bullet=black, label=below:$u^*$] (n8) at (0.65,0.38);
            \coordinate[label=left:$\vdots$] (n9) at (-1.31,0);
            \coordinate[label=right:$\vdots$] (n10) at (1.31,0);
            \draw[-] (n1) -- (n2);
            \draw[-] (n1) -- (n4);
            \draw[-] (n1) -- (n5);
            \draw[-] (n1) -- (n6);
            \draw[-] (n1) -- (n7);
            \draw[-] (n1) -- (n8);
            \draw[-] (n2) -- (n3);
            \draw[-] (n2) -- (n7);
            \draw[-] (n2) -- (n8);
            \draw[-] (n3) -- (n4);
            \draw[-] (n3) -- (n8);
            \draw[-] (n4) -- (n5);
            \draw[-] (n4) -- (n8);
            \draw[-] (n5) -- (n6);
            \draw[-] (n6) -- (n7);
        \end{tikzpicture}
        \caption{A vertex split in a triangulation of a surface}
        \label{fig:steinitzinduction}
    \end{figure}

A \defn{triangulation} of a 2-dimensional manifold $M$ is a pure 2-dimensional simplicial complex $\Sigma$ such that there exists a homeomorphism between $\bigcup\limits_{0\leq k\leq 2}\Sigma^{(k)}$ and $M$.

A \defn{minimal triangulation} of $M$, $\Sigma$, is a triangulation of $M$ for which there is no triangulation of $M$, $\Sigma^{\vee}$, such that $\Sigma$ may be obtained from $\Sigma^{\vee}$ by a vertex split.

Barnette and Edelson showed the following:

\begin{lemma}\cite{barnette1989all}
    Every surface admits finitely many minimal triangulations. 
\end{lemma}

Therefore, if we can show how vertex splitting affects the rigidity of triangulations of surfaces, for each surface, $M$, we only need to consider finitely many triangulations of $M$ to understand the rigidity of all triangulations of $M$.
Indeed, that is what \Cref{lem:vsrig} does.

\begin{lemma}\label{lem:vsrig}
    Let $\Sigma$ and $\Sigma^*$ be two triangulations of a surface.
    Suppose that $\Sigma^*$ is obtained from $\Sigma$ by a vertex split.
    Let $\mathbf{p}\in(\mathbb{R}^2)^n$ and $\mathbf{p}^*\in(\mathbb{R}^2)^{n+1}$ be generic configurations, with $\mathbf{p}^*(v)=\mathbf{p}(v)$, for all $v\in\Sigma^{(0)}$.
    Then $\rank(R(\Sigma^*,\mathbf{p}^*))=\rank(R(\Sigma,\mathbf{p}))+2$.
\end{lemma}

\begin{proof}
    By the rank-nullity theorem, for any matrix $A\in\mathbb{R}^{M\times N}$,
    \begin{equation*}
        \rank(A) = \rank(A^t) = M - \nullity(A^t) = M - \conullity(A),
    \end{equation*}
    where $\conullity(A):=\dim(\coker(A))$.

    Suppose that we split vertex $u$ by deleting 2-simplices $uv_1v_2,\dots,uv_{k-1}v_k$ and adding vertex $u^*$ and 2-simplices $uu^*v_1,uu^*v_k,u^*v_1v_k,\dots,u^*v_{k-1}v_k$, as in \cref{fig:steinitzinduction}.
    Assume that $y_1<\dots<y_t<u<u^*<v_1<\dots<v_k<w_1<\dots<w_{\ell}$.

    Let $\Sigma'$ be the simplicial complex obtained from $\Sigma$ by removing 2-simplices $uv_1v_2,\dots,uv_{k-1}v_k$ (ie. the simplicial complex in the middle panel of \cref{fig:steinitzinduction}).
    The rigidity matrix of $(\Sigma,\mathbf{p})$ has the following form
    \scriptsize
    \begin{equation}\label{eq:splitind1}
        R(\Sigma^*,\mathbf{p}^*)=\begin{bmatrix}
            R(\Sigma',\mathbf{p}) & \mathbf{0} & \mathbf{0} & & \mathbf{A} & & & \mathbf{B} & \\
            \mathbf{0} & \mathbf{n}(u^*v_1,\mathbf{p}^*) & -\mathbf{n}(uv_1,\mathbf{p}^*) & \mathbf{n}(uu^*,\mathbf{p}^*) & \dots & \mathbf{0} & \mathbf{0} & \dots & \mathbf{0} \\
            \mathbf{0} & \mathbf{n}(u^*v_k,\mathbf{p}^*) & -\mathbf{n}(uv_k,\mathbf{p}^*) & \mathbf{0} & \dots & \mathbf{n}(uu^*,\mathbf{p}^*) & \mathbf{0} & \dots & \mathbf{0} \\
            \mathbf{0} & \mathbf{n}(v_1w_{\ell},\mathbf{p}^*) & \mathbf{0} & -\mathbf{n}(uw_{\ell},\mathbf{p}^*) & \dots & \mathbf{0} & \mathbf{0} & \dots & \mathbf{n}(uv_1,\mathbf{p}^*) \\
            \mathbf{0} & \mathbf{n}(v_kw_1,\mathbf{p}^*) & \mathbf{0} & \mathbf{0} & \dots & -\mathbf{n}(uw_1,\mathbf{p}^*) & \mathbf{n}(uv_k,\mathbf{p}^*) & \dots & \mathbf{0} \\
            \mathbf{0} & \mathbf{n}(w_1w_2,\mathbf{p}^*) & \mathbf{0} & \mathbf{0} & \dots & \mathbf{0} & -\mathbf{n}(uw_2,\mathbf{p}^*) & \dots & \mathbf{0} \\
            \vdots & \vdots & \vdots & \vdots & \ddots & \vdots & \vdots & \ddots & \vdots \\ 
            \mathbf{0} & \mathbf{n}(w_{\ell-1}w_{\ell},\mathbf{p}^*) & \mathbf{0} & \mathbf{0} & \dots & \mathbf{0} & \mathbf{0} & \dots & \mathbf{n}(uw_{\ell-1},\mathbf{p}^*) \\
            \mathbf{0} & \mathbf{0} & \mathbf{n}(v_1v_2,\mathbf{p}^*) & -\mathbf{n}(u^*v_2,\mathbf{p}^*) & \dots & \mathbf{0} & -\mathbf{0} & \dots & \mathbf{0} \\
            \vdots & \vdots & \vdots & \vdots & \ddots & \vdots & \vdots & \ddots & \vdots \\ 
            \mathbf{0} & \mathbf{0} & \mathbf{n}(v_{k-1}v_k,\mathbf{p}^*) & \mathbf{0} & \dots & -\mathbf{n}(u^*v_{k-1},\mathbf{p}^*) & \mathbf{0} & \dots & \mathbf{0} \\
        \end{bmatrix}
    \end{equation}
    \normalsize
    where each $\mathbf{0},\mathbf{n}(\tau,\mathbf{p})$ is a $1\times 2$-vector and $\mathbf{A}$ and $\mathbf{B}$ span the column groups of $R(\Sigma,\mathbf{p})$ of size 3 that they lie in the central group of.

    Our aim is to show that each $\omega\in\coker(R(\Sigma,\mathbf{p}))$ uniquely induces $\omega^*\in\coker(R(\Sigma^*,\mathbf{p}^*))$.
    Note that $\coker(R(\Sigma',\mathbf{p}))\subseteq\coker(R(\Sigma,\mathbf{p}))$, so any such $\omega$ must be supported on $S'\sqcup S\subseteq(\Sigma^*)^{(2)}$, where $S'\subseteq(\Sigma')^{(2)}$ and $S\subseteq(\Sigma^*)^{(2)}\setminus(\Sigma')^{(2)}$, with either $|S'|=0$ and $|S|>0$ or $|S'|>0$ and $|S|>0$.

    Firstly, since the columns of $R(\Sigma,\mathbf{p})$ and $R(\Sigma^*,\mathbf{p}^*)$  indexed $y_1,\dots,y_Y$ agree on their non-zero terms $\pi_{S_1}(\omega^*)=\pi_{S_1}(\omega)$, where $S_1 = \{\sigma\in\Sigma^{(2)}:\sigma^{(0)}\cap\{u,v_i,w_j:i\in[k],j\in[\ell]\}=\varnothing\}$.

    The same is true of the columns indexed by $w_1,\dots,w_{\ell}$, so $\pi_{S_2}(\omega^*)=\pi_{S_2}(\omega)$, where $S_2=\{uv_1w_{\ell},uv_kw_1,uw_1w_2,\dots,uw_{\ell-1}w_{\ell}\}$.

    Next, for each $1<i<k-1$,
    \begin{equation}\label{eq:splitind2}
    \begin{split}
        & \mathbf{n}(u^*v_i,\mathbf{p}^*)\omega_{uu^*v_i}^* - \mathbf{n}(u^*v_{i+1},\mathbf{p}^*)\omega_{uu^*v_{i+1}}^* + \sum\limits_{\tau\in\link_{\Sigma}(v_i)\setminus\in\{uv_{i-1},uv_{i+1}\}} \sign(\tau v_i)\omega_{\tau v_i}^*\mathbf{n}(\tau,\mathbf{p}^*) \\
        =& \mathbf{n}(uv_i,\mathbf{p})\omega_{uv_{i-1}v_i} - \mathbf{n}(uv_{i+1},\mathbf{p})\omega_{uuv_iv_{i+1}} + \sum\limits_{\tau\in\link_{\Sigma}(v_i)\setminus\in\{uv_{i-1},uv_{i+1}\}} \sign(\tau v_i)\omega_{\tau v_i}\mathbf{n}(\tau,\mathbf{p}) \\
        =& \mathbf{0},
    \end{split}
    \end{equation}
    the summands of rightmost terms in each line of \cref{eq:splitind2} are pairwise equal, so we must have that
    \begin{equation}\label{eq:splitind3}
        \mathbf{n}(u^*v_i,\mathbf{p}^*)\omega_{uu^*v_i}^* - \mathbf{n}(u^*v_{i+1},\mathbf{p}^*)\omega_{uu^*v_{i+1}}^* = \mathbf{n}(uv_i,\mathbf{p})\omega_{uv_{i-1}v_i} - \mathbf{n}(uv_{i+1},\mathbf{p})\omega_{uuv_iv_{i+1}}.
    \end{equation}
    Recall that $\mathbf{n}(\tau,\mathbf{p})\in\mathbb{R}^2$, so \cref{eq:splitind3} is two equations, one in each coordinate, which, by the genericity of $\mathbf{p}^*$, are independent.
    These two independent equations uniquely define $\omega_{uu^*v_i}^*$ and $\omega_{uu^*v_{i+1}}^*$ as
    \begin{equation}
        \begin{bmatrix}
            \omega_{uu^*v_i}^* \\ \omega_{uu^*v_{i+1}}^*
        \end{bmatrix} = \begin{bmatrix}
            -\omega_{uv_iv_{i+1}} \\ \omega_{uv_iv_{i+1}}
        \end{bmatrix}.
    \end{equation}

    This leaves four as-of-yet undefined entries of $\omega^*$: $\omega_{uu^*v_1}^*,\omega_{uu^*v_k}^*,\omega_{uv_1,w_{\ell}}^*$ and $\omega_{uv_kw_1}^*$.
    Considering the product of $\omega^*$ (resp. $\omega$) with column group $v_1$ of $R(\Sigma^*,\mathbf{p}^*)$ (resp. $R(\Sigma,\mathbf{p})$) gives us the following
    \begin{equation}\label{eq:splitind5}
    \begin{split}
        & \omega_{uu^*v_1}^*\mathbf{n}(uu^*,\mathbf{p}^*) - \omega_{u^*v_1v_2}^*\mathbf{n}(u^*v_2,\mathbf{p}^*) + \sum\limits_{\substack{\tau\in\link_{\Sigma}(v_1) \\ \tau\not\in\{uu^*v_1,u^*v_1v_2\}}}\sign(v_1,\tau v_1)\omega_{v_1\tau}^*\mathbf{n}(\tau,\mathbf{p}^*) \\
        =& - \omega_{uv_1v_2}\mathbf{n}(uv_2,\mathbf{p}) + \sum\limits_{\substack{\tau\in\link_{\Sigma}(v_1) \\ \tau\not\in\{uv_1v_2\}}}\sign(v_1,\tau v_1)\omega_{v_1\tau}\mathbf{n}(\tau,\mathbf{p}) \\
        =& \mathbf{0},
    \end{split}
    \end{equation}
    which, by the pairwise equality of the summands of the rightmost terms of both the first two lines of \cref{eq:splitind5}, we get
    \begin{equation}\label{eq:splitind6}
        \omega_{uu^*v_1}^*\mathbf{n}(uu^*,\mathbf{p}^*) - \omega_{u^*v_1v_2}^*\mathbf{n}(u^*v_2,\mathbf{p}^*) = -\omega_{uv_1v_2}\mathbf{n}(uv_2,\mathbf{p}),
    \end{equation}
    and similarly for column group $v_k$:
    \begin{equation}\label{eq:splitind7}
        \omega_{uu^*v_k}^*\mathbf{n}(uu^*,\mathbf{p}^*) + \omega_{u^*v_{k-1}v_k}^*\mathbf{n}(u^*v_{k-1},\mathbf{p}^*) = \omega_{uv_{k-1}v_k}\mathbf{n}(uv_{k-1},\mathbf{p}).
    \end{equation}
    Each of \cref{eq:splitind6} and \cref{eq:splitind7} denotes two linearly independent linear equations in two variables.
    Therefore, between them, they uniquely define $\omega_{uu^*v_1}^*, \omega_{uu^*v_k}^*,\omega_{u^*v_1v_2}^*,\omega_{u^*v_{k-1}k}^*$ as follows:
    \begin{equation}\label{eq:splitind8}
        \begin{bmatrix}
            \omega_{uu^*v_1}^* \\ \omega_{uu^*v_k}^* \\ \omega_{u^*v_1v_2}^* \\ \omega_{u^*v_{k-1}v_k}^*
        \end{bmatrix}
        = \begin{bmatrix}
            \omega_{uv_1v_2} \\ \omega_{uv_{k-1}v_k} \\ -\omega_{uv_1v_2} \\ \omega_{uv_{k-1}v_k}
        \end{bmatrix}.
    \end{equation}
    Therefore $\Span\{\omega^*:\omega\in\coker(R(\Sigma,\mathbf{p}))\}\subseteq\coker(R(\Sigma^*,\mathbf{p}^*))$.

    It remains to show equality.
    
    Suppose that $\omega^*\in\coker(R(\Sigma^*,\mathbf{p}^*))$, but $\pi_{S_1}(\omega^*)=0$.

    Let $1<i<k-2$, then if $\omega_{u^*v_iv_{i+1}}^*\neq 0$, we have that
    \begin{equation}\label{eq:splitind9}
        \omega_{u^*v_iv_{i+1}}^*\mathbf{n}(u^*v_i,\mathbf{p}^*) - \omega_{u^*v_{i+1}v_{i+2}}^*\mathbf{n}(u^*v_{i+2},\mathbf{p}^*) = 0.
    \end{equation}
    These two equations contradict the genericity of $\mathbf{p}$, so $\omega_{u^*v_iv_{i+1}}^*=0$, for all $1<i<k-2$.
    By a similar argument, $\omega_{u^*w_jw_{j+1}}^*=0$, for all $1<j<\ell-2$.

    This leaves four as-of-yet undetermined entries of $\omega^*$: $\omega_{uu^*v_1}^*,\omega_{uu^*v_k}^*,\omega_{uv_1w_{\ell}}^*,\omega_{uv_kw_1}^*$.
    The first two may be determined as follows:
    \begin{equation}
        -\omega_{uu^*v_1}^*\mathbf{n}(uv_1,\mathbf{p}^*) + \sum\limits_{i=1}^k-1\omega_{u^*v_iv_{i+1}}^*\mathbf{n}(v_iv_{i+1},\mathbf{p}^*)
        = -\omega_{uu^*v_1}^*\mathbf{n}(uv_1,\mathbf{p}^*)
        = 0,
    \end{equation}
    and thus $\omega_{uu^*v_1}^*=0$, similarly, $\omega_{uu^*v_k}^*=0$.
    Finally,
    \begin{equation}
        \omega_{uv_1w_{\ell}}^*\mathbf{n}(v_1w_{\ell},\mathbf{p}^*) + \omega_{uv_kw_1}^*\mathbf{n}(v_kw_1,\mathbf{p}^*) = 0
    \end{equation}
    forces $\omega_{uv_1w_{\ell}}^*=\omega_{uv_kw_1}^*=0$.
    Therefore, if $\pi_{S_1}(\omega^*)=0$, then $\omega^*=0$.
\end{proof}

Prior to Barnette and Edelson, Steinitz showed the following:

\begin{lemma}\label{lem:steinitz}\cite{steinitz2013vorlesungen}
    The sphere $\mathbb{S}^2$ admits a single minimal triangulation: the tetrahedron $K_4^2$.
\end{lemma}

Combining \cref{prop:compgr}, \Cref{lem:vsrig} and \Cref{lem:steinitz}, we see that all triangulations of $\mathbb{S}^2$ are rigid in $\mathbb{R}^2$, formalised in the following proposition:

\begin{proposition}\label{prop:rigs2}
    Let $\Sigma$ be a triangulation of $\mathbb{S}^2$.
    Then $\Sigma$ is rigid in $\mathbb{R}^2$.
\end{proposition}

Bulavka et al. \cite{bulavka2022volume} show that the Klein bottle and the torus are also volume rigid by considering their minimal triangulations (of which there are 21 and 29 respectively).

Can we go further to say that generic frameworks of some (or all) triangulations of $\mathbb{S}^2$ are globally rigid in $\mathbb{R}^2$?

A \defn{stacked tetrahedron} is a triangulation of $\mathbb{S}^2$ obtained by performing repeated vertex splits to $K_4^2$ where one 2-simplex is removed and three 2-simplices are added in its place.

\begin{proposition}\label{lem:stackedtet}
    Let $\Sigma$ be a stacked tetrahedron.
    Any generic framework of $\Sigma$ in $\mathbb{R}^2$ is globally rigid. 
\end{proposition}

\begin{proof}
    We will proceed by induction.
    We showed in \cref{prop:compgr} that any generic framework of $K_4^2$ is globally rigid in $\mathbb{R}^2$.
    Let $K_4^2,p$ be such a framework.

    Now let $\Sigma$ and $\Sigma^*$ be stacked tetrahedra on $n$ and $n+1$ vertices respectively, with $\Sigma^*$ obtained from $\Sigma$ by splitting vertex $n$ to obtain vertex $n+1$, deleting 2-simplex $12n$ and adding 2-simplices $12(n+1), 1n(n+1), 2n(n+1)$.

    Assume that any generic framework of $\Sigma$ in $\mathbb{R}^2$ is globally rigid. Let $(\Sigma^*,\mathbf{p}^*)$ be a generic framework of $\Sigma^*$ in $\mathbb{R}^2$.
    Pin $(\Sigma^*,\mathbf{p}^*)$ so that $\mathbf{p}^*(1)=(0,0)$, $\mathbf{p}^*(2)=(1,0)$ and $\mathbf{p}^*(3)=(0,1)$ to obtain $\overline{\mathbf{p}^*}$.
    Let $(\Sigma,\overline{\mathbf{p}})$ be the pinned framework obtained by projecting $\overline{\mathbf{p}^*}$ onto its first $2n$ coordinates.
    Although $(\Sigma,\overline{\mathbf{p}})$ is not generic, it is related to a generic framework by the composition of a special affine transformation of $\mathbb{R}^2$ and scaling in one coordinate direction, transformations which do not affect its (global) rigidity theoretic properties \cite[p. 56]{southgate2024volume}.

    Now, suppose that $(\Sigma^*,\mathbf{p}^*)$ and $(\Sigma^*,\mathbf{q})$ are equivalent generic frameworks, then the line arrangement defining the position of $\mathbf{q}(n+1)$ in $(\Sigma^*,\mathbf{q})$ is non-degenerately transformed to an arrangement uniquely defining $\overline{\mathbf{q}}(n+1)$ in terms of $\{\overline{\mathbf{p}^*}(i),\overline{\mathbf{q}}(j):i\in [n+1], j\in [n]\}$.
    By the global rigidity of $(\Sigma,\mathbf{p})$, the position of each $\overline{\mathbf{q}}(j)$ is uniquely defined as being equal to $\overline{\mathbf{p}^*}(j)$.
    Therefore, $\overline{\mathbf{q}}(n+1)=\overline{\mathbf{p}}^*(n+1)$, and so $(\Sigma^*,\overline{\mathbf{q}})=(\Sigma^*,\overline{\mathbf{p}^*})$.
\end{proof}

In other words, for any $n\geq 4$, there exists a triangulation of $\mathbb{S}^2$, every generic framework of which is globally rigid in $\mathbb{R}^2$.

\section{Bipyramids}

A \defn{bipyramid} on $n$ vertices, denoted $B_{n-2}$, is a 2-dimensional simplicial complex defined by its maximal simplices
\begin{equation*}
    B_{n-2}^{(2)} = \{123, 23n, \dots, 1(n-2)(n-1), (n-2)(n-1)n, 12(n-1), 2(n-1)n\}.
\end{equation*}

For brevity, we will refer to vertices 1 and $n$ of $B_{n-2}$ as its \textit{south} and \textit{north pole}, and as vertices $2,\dots, n-1$ as lying on its equator, or being \textit{equatorial}.
The reason for this terminology can be seen in \cref{fig:bipyramid}.

\begin{figure}
        \centering
        \begin{tikzpicture}
            \filldraw[fill=yellow!80!black,line width=0.5pt,fill opacity=0.6] (0.2,0) -- (-1,1) -- (1,1) -- (0.2,0);
            \filldraw[fill=yellow!80!black,line width=0.5pt,fill opacity=0.6] (0.2,0) -- (-1,1) -- (-0.4,1.2) -- (0.2,0);
            \filldraw[fill=yellow!80!black,line width=0.5pt,fill opacity=0.6] (0.2,0) -- (1,1) -- (1.4,1.2) -- (0.2,0);
            \filldraw[fill=yellow!80!black,line width=0.5pt,fill opacity=0.6] (0.2,0) -- (-0.4,1.2) -- (1.4,1.2) -- (0.2,0);
            \filldraw[fill=yellow!80!black,line width=0.5pt,fill opacity=0.6] (0.2,2.2) -- (-1,1) -- (1,1) -- (0.2,2.2);
            \filldraw[fill=yellow!80!black,line width=0.5pt,fill opacity=0.6] (0.2,2.2) -- (-1,1) -- (-0.4,1.2) -- (0.2,2.2);
            \filldraw[fill=yellow!80!black,line width=0.5pt,fill opacity=0.6] (0.2,2.2) -- (1,1) -- (1.4,1.2) -- (0.2,2.2);
            \filldraw[fill=yellow!80!black,line width=0.5pt,fill opacity=0.6] (0.2,2.2) -- (-0.4,1.2) -- (1.4,1.2) -- (0.2,2.2);
            \coordinate[Bullet=black, label=below:1] (n1) at (0.2,0);
            \coordinate[Bullet=black, label=below:2] (n2) at (-1,1);
            \coordinate[Bullet=black, label=below:3] (n3) at (1,1);
            \coordinate[Bullet=black, label=below:{4}] (n4) at (-0.4,1.2);
            \coordinate[Bullet=black, label=below:{5}] (n5) at (1.4,1.2);
            \coordinate[Bullet=black, label=above:{6}] (n6) at (0.2,2.2);
        \end{tikzpicture}
        \caption{The octahedron, or $B_4$}
        \label{fig:bipyramid}
\end{figure}

For brevity, let an \defn{$m$-vertex split} in a triangulation of a surface be one that deletes $m-2$ 2-simplices and adds $m$ back.

\begin{lemma}\label{lem:bipsplit}
    There is a sequence of one 3-vertex split and $n-5$ 4-vertex splits from $K_4^2$ to $B_{n-2}$.
\end{lemma}

\begin{proof}
    Begin with $K_4^2$, there is a single 3-vertex split we can perform, up to isomorphism.
    Without loss of generality, split vertex 4, deleting 2-simplex 234 and creating a new vertex 5 and new 2-simplices 235, 245 and 345.
    This yields the bipyramid $B_3$.

    At this point, there are two vertex splits that we may perform, up to isomorphism: a 3-vertex split at vertex 5, which we will not consider, and a 4-vertex split at vertex 4, deleting 2-simplices 124 and 245 and creating a new vertex $4^*$ and new 2-simplices $124^*, 144^*, 24^*5$ and $44^*5$.
    We then relabel vertices 5 and $4^*$ as 6 and 5 respectively to obtain $B_4$.

    Repeating 4-vertex splits along the equator in this manner yields bipyramids on successively larger vertex sets, increasing in increments of one.
\end{proof}

We now restate and prove \cref{thm:main}.

\begin{theorem}\label{thm:bipcount}
    For each $n\geq 5$, $c(B_{n-2})\leq n-4$.
\end{theorem}

\begin{proof}
    To prove this theorem, we will take a generic framework $(B_{n-2},\mathbf{p})$, pin it to obtain $(B_{n-2},\overline{\mathbf{p}})$ and show that a pinned framework $(B_{n-2},\overline{\mathbf{q}})$ is equivalent to $(B_{n-2},\overline{\mathbf{p}})$ if and only if it satisfies a polynomial equation with coefficients in $\mathbb{Q}[\overline{\mathbf{p}}]$ of degree $n-4$.

    Pin $(B_{n-2},\mathbf{p})$ at 2-simplex $123$ so that $\overline{\mathbf{p}}(1)=(0,0)$, $\overline{\mathbf{p}}(2)=(1,0)$ and $\overline{\mathbf{p}}(3)=(0,1)$.
    Note again that, although this pinned framework is not generic, scaling by a generic quantity in the $y$-direction then performing a suitable special affine transformation of $\mathbb{R}^2$ yields a generic framework without changing any of its rigidity-theoretic properties.

    Suppose that $(B_{n-2},\overline{\mathbf{p}})$ is equivalent to $(B_{n-2},\overline{\mathbf{q}})$.
    Then, $\det(C(123,\overline{\mathbf{p}}))=\det(C(123,\overline{\mathbf{q}})$ (immediately, as they are both pinned to the same coordinates).
    The following three equalities define three lines on which the corresponding points of $(B_{n-2},\overline{\mathbf{q}})$ must lie:
    \begin{equation}
    \begin{split}
        \det(C(12(n-1),\overline{\mathbf{p}})) &= \det(C(12(n-1),\overline{\mathbf{p}})),\\
        \det(C(134,\overline{\mathbf{p}})) &= \det(C(134,\overline{\mathbf{p}})),\\
        \det(C(23n,\overline{\mathbf{p}})) &= \det(C(234,\overline{\mathbf{p}}));
    \end{split}
    \end{equation}
    respectively:
    \begin{equation}\label{eq:bip1}
        \overline{\mathbf{q}}(n-1) = \begin{bmatrix} \overline{\mathbf{p}}(n-1)_1 + s \\ \overline{\mathbf{p}}(n-1)_2 \end{bmatrix}\text{, }
        \overline{\mathbf{q}}(4) = \begin{bmatrix} \overline{\mathbf{p}}(4)_1 \\ \overline{\mathbf{p}}(4)_2 + t\end{bmatrix}\text{, }
        \overline{\mathbf{q}}(n) = \begin{bmatrix} \overline{\mathbf{p}}(n)_1 + r \\ \overline{\mathbf{p}}(n)_2 - r \end{bmatrix},
    \end{equation}
    with $r,s,t$ varying in $\mathbb{R}$.

    Next, plugging in the values from \cref{eq:bip1} into the equations $\det(C(2(n-1)n,\overline{\mathbf{p}}))=\det(C(2(n-1)n,\overline{\mathbf{q}}))$ and $\det(C(34n,\overline{\mathbf{p}}))=\det(C(34n,\overline{\mathbf{q}}))$, we obtain
    \begin{equation}\label{eq:bip2}
    \begin{split}
        r&=\frac{\overline{\mathbf{p}}(n)_1t}{1-\overline{\mathbf{p}}(4)_1-\overline{\mathbf{p}}(4)_2-t}=\frac{\overline{\mathbf{p}}(n)_2s}{1-\overline{\mathbf{p}}(n-1)_1-\overline{\mathbf{p}}(n-1)_2-2} \\
        \implies&s=\frac{(\overline{\mathbf{p}}(n-1)_1+\overline{\mathbf{p}}(n-1)_2-1)\overline{\mathbf{p}}(n)_1t}{(\overline{\mathbf{p}}(4)_1+\overline{\mathbf{p}}(4)_2-1)\overline{\mathbf{p}}(n)_2+(\overline{\mathbf{p}}(n)_1+\overline{\mathbf{p}}(n)_2)t},
    \end{split}
    \end{equation}
    thus putting all our equations in terms of one variable: $t$.

    Next, we define the positions of the equatorial vertices of $(B_{n-2},\overline{\mathbf{q}})$ in terms of $\overline{\mathbf{p}}$ and $t$.
    For each $4\leq i\leq n-1$, the position of $\overline{\mathbf{q}}(i)$, is completely determined by $\overline{\mathbf{q}}(1)=(0,0)$, $\overline{\mathbf{q}}(i-1)$ and $\overline{\mathbf{q}}(n)=(\overline{\mathbf{p}}(n)_1+r,\overline{\mathbf{p}}(n)_2-r)$.
    Therefore,
    \begin{equation}\label{eq:bip3}
    \begin{split}
        \overline{\mathbf{q}}(i)_j=&\frac{\left(\begin{vmatrix} \overline{\mathbf{q}}(i-1)_1 & \overline{\mathbf{q}}(n)_1 \\ \overline{\mathbf{q}}(i-1)_2 & \overline{\mathbf{q}}(n)_2 \end{vmatrix} - \begin{vmatrix} \overline{\mathbf{p}}(i-1)_1 & \overline{\mathbf{p}}(n)_1 \\ \overline{\mathbf{p}}(i-1)_2 & \overline{\mathbf{p}}(n)_2 \end{vmatrix} \begin{vmatrix} \overline{\mathbf{p}}(i)_1 & \overline{\mathbf{p}}(n)_1 \\ \overline{\mathbf{p}}(i)_2 & \overline{\mathbf{p}}(n)_2 \end{vmatrix} \right)\overline{\mathbf{q}}(i-1)_j}{\begin{vmatrix} \overline{\mathbf{q}}(i-1)_1 & \overline{\mathbf{q}}(n)_1 \\ \overline{\mathbf{q}}(i-1)_2 & \overline{\mathbf{q}}(n)_2 \end{vmatrix}} \\
        &+ \frac{\begin{vmatrix} \overline{\mathbf{p}}(i-1)_1 & \overline{\mathbf{p}}(n)_1 \\ \overline{\mathbf{p}}(i-1)_2 & \overline{\mathbf{p}}(n)_2 \end{vmatrix}\overline{\mathbf{q}}(n)_j}{\begin{vmatrix} \overline{\mathbf{q}}(i-1)_1 & \overline{\mathbf{q}}(n)_1 \\ \overline{\mathbf{q}}(i-1)_2 & \overline{\mathbf{q}}(n)_2 \end{vmatrix}},
    \end{split}
    \end{equation}
    for each $j\in[2]$.

    In order to simplify \cref{eq:bip3}, notice the following identities:
    \begin{equation}\label{eq:detind1}
        \begin{vmatrix} \overline{\mathbf{p}}(i)_1 & \overline{\mathbf{p}}(n)_1 \\ \overline{\mathbf{p}}(i-1)_2 & \overline{\mathbf{p}}(n)_2 \end{vmatrix} - \begin{vmatrix} \overline{\mathbf{p}}(i-1)_1 & \overline{\mathbf{p}}(n)_1 \\ \overline{\mathbf{p}}(i-1)_2 & \overline{\mathbf{p}}(n)_2 \end{vmatrix} = \det(C(1in,\overline{\mathbf{p}})) = \det(C(1(i-1)n,\overline{\mathbf{p}}))
    \end{equation}
    and
    \begin{equation}\label{eq:detind2}
        \det(C(1(i-1)i,\overline{\mathbf{p}})) - \det(C((i-1)in,\overline{\mathbf{p}})) = \det(C(1(i-1)n,\overline{\mathbf{p}})) - \det(C(1in,\overline{\mathbf{p}})).
    \end{equation}
    The right hand side of \cref{eq:detind2} is constant across equivalent pinned frameworks, hence
    \begin{equation}\label{eq:bip4}
        \begin{vmatrix} \overline{\mathbf{q}}(i)_1 & \overline{\mathbf{q}}(n)_1 \\ \overline{\mathbf{q}}(i)_2 & \overline{\mathbf{q}}(n)_2 \end{vmatrix} - \begin{vmatrix} \overline{\mathbf{q}}(i-1)_1 & \overline{\mathbf{q}}(n)_1 \\ \overline{\mathbf{q}}(i-1)_2 & \overline{\mathbf{q}}(n)_2 \end{vmatrix} = \begin{vmatrix} \overline{\mathbf{p}}(i)_1 & \overline{\mathbf{p}}(n)_1 \\ \overline{\mathbf{p}}(i)_2 & \overline{\mathbf{p}}(n)_2 \end{vmatrix} - \begin{vmatrix} \overline{\mathbf{p}}(i-1)_1 & \overline{\mathbf{p}}(n)_1 \\ \overline{\mathbf{p}}(i-1)_2 & \overline{\mathbf{p}}(n)_2 \end{vmatrix}
    \end{equation},
    for each $4\leq i\leq n-1$.
    Next, after applying \cref{eq:bip4} sufficiently many times, we obtain
    \begin{equation}
    \begin{split}
        \begin{vmatrix}
            \overline{\mathbf{q}}(i-1)_1 & \overline{\mathbf{q}}(n)_1 \\ \overline{\mathbf{q}}(i-1)_2 & \overline{\mathbf{q}}(n)_2
        \end{vmatrix} &= \begin{vmatrix} \overline{\mathbf{q}}(3)_1 & \overline{\mathbf{q}}(n)_1 \\ \overline{\mathbf{q}}(3)_2 & \overline{\mathbf{q}}(n)_2 \end{vmatrix} - \begin{vmatrix} \overline{\mathbf{p}}(3)_1 & \overline{\mathbf{p}}(n)_1 \\ \overline{\mathbf{p}}(3)_2 & \overline{\mathbf{p}}(n)_2 \end{vmatrix} + \begin{vmatrix} \overline{\mathbf{p}}(i-1)_1 & \overline{\mathbf{p}}(n)_1 \\ \overline{\mathbf{p}}(i-1)_2 & \overline{\mathbf{p}}(n)_2 \end{vmatrix} \\
        &= \begin{vmatrix} 0 & \overline{\mathbf{p}}(n)_1 + r \\ 1 & \overline{\mathbf{p}}(n)_2 - r \end{vmatrix} - \begin{vmatrix} 0 & \overline{\mathbf{p}}(n)_1 \\ 1 & \overline{\mathbf{p}}(n)_2 \end{vmatrix} + \begin{vmatrix} \overline{\mathbf{p}}(i-1)_1 & \overline{\mathbf{p}}(n)_1 \\ \overline{\mathbf{p}}(i-1)_2 & \overline{\mathbf{p}}(n)_2 \end{vmatrix} \\
        &= \begin{vmatrix}
            \overline{\mathbf{p}}(i-1)_1 & \overline{\mathbf{p}}(n)_1 \\ \overline{\mathbf{p}}(i-1)_2 & \overline{\mathbf{p}}(n)_2
        \end{vmatrix} - r.
    \end{split} 
    \end{equation}
    Plugging this into our formula from \cref{eq:bip3}, we obtain
    \begin{equation}\label{eq:bip5}
        \overline{\mathbf{q}}(i)_j = \frac{\left(\begin{vmatrix} \overline{\mathbf{p}}(i)_1 & \overline{\mathbf{p}}(n)_1 \\ \overline{\mathbf{p}}(i)_2 & \overline{\mathbf{p}}(n)_2 \end{vmatrix} - r\right)\overline{\mathbf{q}}(i-1)_j + \begin{vmatrix} \overline{\mathbf{p}}(i-1)_1 & \overline{\mathbf{p}}(i)_1 \\ \overline{\mathbf{p}}(i-1)_2 & \overline{\mathbf{p}}(i)_2 \end{vmatrix}(\overline{\mathbf{p}}(n)_j - (-1)^jr)}{\begin{vmatrix} \overline{\mathbf{p}}(i-1)_1 & \overline{\mathbf{p}}(i)_1 \\ \overline{\mathbf{p}}(i-1)_2 & \overline{\mathbf{p}}(i)_2 \end{vmatrix} - r},
    \end{equation}
    for each $j\in[2]$ and $4\leq i \leq n-1$.

    In order for $(B_{n-2},\overline{\mathbf{q}})$ as we are defining it to be a well-defined framework, the position of $\overline{\mathbf{q}}(n-1)$ obtained from \cref{eq:bip5} must line up with that from \cref{eq:bip1}.
    This equality is obtained when the following two equations are satisfied:
    \begin{equation}\label{eq:bip6}
    \begin{split}
        &\left(\left( \begin{vmatrix} \overline{\mathbf{p}}(n-1)_1 & \overline{\mathbf{p}}(n)_1 \\ \overline{\mathbf{p}}(n-1)_2 & \overline{\mathbf{p}}(n)_2 \end{vmatrix} (1-\overline{\mathbf{p}}(4)_1 - \overline{\mathbf{p}}(4)_2 -t) - \overline{\mathbf{p}}(n)_1t \right)\overline{\mathbf{q}}(n-2)_1\right.\\
        &+\left. \begin{vmatrix} \overline{\mathbf{p}}(n-2)_1 & \overline{\mathbf{p}}(n-1)_1 \\ \overline{\mathbf{p}}(n-2)_2 & \overline{\mathbf{p}}(n-2)_2 \end{vmatrix} (1-\overline{\mathbf{p}}(4)_1-\overline{\mathbf{p}}(4)_2)\overline{\mathbf{p}}(n)_1 \right) \\
        &((\overline{\mathbf{p}}(4)_1+\overline{\mathbf{p}}(4)_2 - 1)\overline{\mathbf{p}}(n)_2 - (\overline{\mathbf{p}}(n)_1+\overline{\mathbf{p}}(n)_2)t) \\
        =& \left( \begin{vmatrix} \overline{\mathbf{p}}(n-2)_1 & \overline{\mathbf{p}}(n-1)_1 \\ \overline{\mathbf{p}}(n-2)_2 & \overline{\mathbf{p}}(n-1)_2 \end{vmatrix} (1-\overline{\mathbf{p}}(4)_1-\overline{\mathbf{p}}(4)_2-t) - \overline{\mathbf{p}}(n)_1t \right) \\
        &\left( (\overline{\mathbf{p}}(4)_1+\overline{\mathbf{p}}(4)_2-1)\overline{\mathbf{p}}(n-1)_1\overline{\mathbf{p}}(n)_2 + \left(\begin{vmatrix} \overline{\mathbf{p}}(n-1)_1 & \overline{\mathbf{p}}(n)_1 \\ \overline{\mathbf{p}}(n-1)_2 & \overline{\mathbf{p}}(n)_2 \end{vmatrix}\right)t \right)
    \end{split}
    \end{equation}
    and
    \begin{equation}\label{eq:bip7}
    \begin{split}
        &\left(\left( \begin{vmatrix} \overline{\mathbf{p}}(n-1)_1 & \overline{\mathbf{p}}(n)_1 \\ \overline{\mathbf{p}}(n-1)_2 & \overline{\mathbf{p}}(n)_2 \end{vmatrix} (1-\overline{\mathbf{p}}(4)_1 - \overline{\mathbf{p}}(4)_2 -t) - \overline{\mathbf{p}}(n)_1t \right)\overline{\mathbf{q}}(n-2)_2\right.\\
        &+\left. \begin{vmatrix} \overline{\mathbf{p}}(n-2)_1 & \overline{\mathbf{p}}(n-1)_1 \\ \overline{\mathbf{p}}(n-2)_2 & \overline{\mathbf{p}}(n-2)_2 \end{vmatrix} ((1-\overline{\mathbf{p}}(4)_1-\overline{\mathbf{p}}(4)_2)\overline{\mathbf{p}}(n)_2 - \overline{\mathbf{p}}(n)_1t) \right) \\
        =& \left( \begin{vmatrix} \overline{\mathbf{p}}(n-2)_1 & \overline{\mathbf{p}}(n-1)_1 \\ \overline{\mathbf{p}}(n-2)_2 & \overline{\mathbf{p}}(n-1)_2 \end{vmatrix} (1-\overline{\mathbf{p}}(4)_1-\overline{\mathbf{p}}(4)_2-t) - \overline{\mathbf{p}}(n)_1t \right) \overline{\mathbf{p}}(n-1)_2.
    \end{split}
    \end{equation}
    Notice, however, that it suffices to solve just \cref{eq:bip7}.
    Indeed, if $\tau$ solves \cref{eq:bip7}, so $\overline{\mathbf{q}}(\tau)(n-1)_2=\overline{\mathbf{p}}(n-1)_2$, but not \cref{eq:bip6}, so $\overline{\mathbf{q}}(\tau)(n-1)_1\neq\overline{\mathbf{p}}(n-1)_1+\sigma$, for some $\sigma\neq s$, then comparing the value of $r$ yielded by $\det(C(2(n-1)n,\overline{\mathbf{p}}))=\det(C(2(n-1)n,\overline{\mathbf{q}}(\tau)))$ gives us
    \begin{equation*}
        s(\overline{\mathbf{p}}(n-1)_1+\overline{\mathbf{p}}(n-1)_2-1)=\sigma(\overline{\mathbf{p}}(n-1)_1+\overline{\mathbf{p}}(n-1)_2),
    \end{equation*}
    hence $s=\sigma$, a contradiction.

    Subtracting the right hand side from both sides of \cref{eq:bip7} yields an equation of the form $\text{LHS}=0$, then multiplying through by the denominator of $\overline{\mathbf{q}}(n-2)_2$ increases the degree of terms not containing $\overline{\mathbf{q}}(n-2)_2$ by 1 (with respect to the variable $t$).
    This yields an equation of the form $f(\overline{\mathbf{p}})(t)=0$, with
    \begin{equation*}
        \deg(f(\mathbf{p}))=\begin{cases}
            \max\{\deg(\overline{\mathbf{q}}(n-2)_2)+1,2\},\text{ if }n\geq 6,\\
            \deg(\overline{\mathbf{q}}(n-2)_2+1,\text{ if }n=5.
        \end{cases}
    \end{equation*}
    By considering \cref{eq:bip5}, we notice that the degree of the numerator of $\overline{\mathbf{q}}(i)_2$ has a higher degree than its denominator, and so, multiplying $f(\mathbf{p})(t)=0$ through by the denominator of $\overline{\mathbf{q}}(n-2)_2$ does not change the degree of $f(\mathbf{p})$, hence the numerator's degree increases by 1 for each equatorial vertex, with an initial value of 0 at $\overline{\mathbf{q}}(3)_2$.
    Therefore $\deg(f(\mathbf{p}))=n-4$.

    This completes the proof subject to one final claims:
        If $(B_{n-2},\mathbf{p})$ and $(B_{n-2},\mathbf{q}_1)$ and $(B_{n-2},\mathbf{q}_2)$ are equivalent and if $(B_{n-2},\overline{\mathbf{q}_1})\ne\mathbf{q}((B_{n-2},\overline{\mathbf{q}_2})$, then $t_1=\overline{\mathbf{q}_1}(4)_2-\overline{\mathbf{p}}(4)_2\neq\overline{\mathbf{q}_2}(4)_2-\overline{\mathbf{p}}(4)_2=t_2$.

        Indeed, if, to the contrary, $t_1=t_2$, then $r_1=r_2$ and $s_1=s_2$ (with $r_1,r_2,s_1,s_2$ analogously defined).
        Therefore each equatorial vertex of $(B_{n-2},\overline{\mathbf{q}_2})$ has the same position of its corresponding vertex in $(B_{n-2},\overline{\mathbf{q}_1})$, hence $(B_{n-2},\overline{\mathbf{q}_1})=(B_{n-2},\overline{\mathbf{q}_2})$.
\end{proof}

\Cref{fig:bipcount} outlines the process of the proof of \cref{thm:bipcount} pictorially.

\begin{figure}
        \centering
        \begin{tikzpicture}
            \filldraw[fill=yellow!80!black,line width=0.5pt,fill opacity=0.6] (0,0) -- (3,0) -- (0,3) -- (0,0);
            \filldraw[fill=yellow!80!black,line width=0.5pt,fill opacity=0.6] (0,0) -- (3,0) -- (1.4,0.5) -- (0,0);
            \filldraw[fill=yellow!80!black,line width=0.5pt,fill opacity=0.6] (0,0) -- (0,3) -- (0.5,1.4) -- (0,0);
            \filldraw[fill=yellow!80!black,line width=0.5pt,fill opacity=0.6] (3,0) -- (0,3) -- (1.2,1.2) -- (3,0);
            \coordinate[Bullet=black, label=below:1] (n1) at (0,0);
            \coordinate[Bullet=black, label=below:2] (n2) at (3,0);
            \coordinate[Bullet=black, label=left:3] (n3) at (0,3);
            \coordinate[Bullet=black, label=left:{4}] (n4) at (0.5,1.4);
            \coordinate[Bullet=black, label=below:{$n-1$}] (n5) at (1.4,0.5);
            \coordinate[Bullet=black, label=right:{$n$}] (n6) at (1.2,1.2);
            \draw[dashed,red] (0,0.5) -- (3,0.5);
            \draw[dashed,red] (0.5,0) -- (0.5,3);
            \draw[dashed,red] (0,2.4) -- (2.4,0);
        \end{tikzpicture}
        \begin{tikzpicture}
            \filldraw[fill=yellow!80!black,line width=0.5pt,fill opacity=0.6] (0,0) -- (3,0) -- (0,3) -- (0,0);
            \filldraw[fill=yellow!80!black,line width=0.5pt,fill opacity=0.6] (0,0) -- (3,0) -- (1.4,0.5) -- (0,0);
            \filldraw[fill=yellow!80!black,line width=0.5pt,fill opacity=0.6] (0,0) -- (0,3) -- (0.5,1.4) -- (0,0);
            \filldraw[fill=yellow!80!black,line width=0.5pt,fill opacity=0.6] (3,0) -- (0,3) -- (1.2,1.2) -- (3,0);
            \filldraw[fill=yellow!80!black,line width=0.5pt,fill opacity=0.6] (3,0) -- (1.4,0.5) -- (1.2,1.2) -- (3,0);
            \filldraw[fill=yellow!80!black,line width=0.5pt,fill opacity=0.6] (0,3) -- (0.5,1.4) -- (1.2,1.2) -- (0,3);
            \coordinate[Bullet=black, label=below:1] (n1) at (0,0);
            \coordinate[Bullet=black, label=below:2] (n2) at (3,0);
            \coordinate[Bullet=black, label=left:3] (n3) at (0,3);
            \coordinate[Bullet=black, label=left:{4}] (n4) at (0.5,1.4);
            \coordinate[Bullet=black, label=below:{$n-1$}] (n5) at (1.4,0.5);
            \coordinate[Bullet=black, label=right:{$n$}] (n6) at (1.2,1.2);
            \draw[->,red] (1.4,0.5) -- (2,0.5);
            \draw[->,red] (0.5,1.4) -- (0.5,0.8);
            \draw[->,red] (1.2,1.2) -- (0.8,1.6);
        \end{tikzpicture}
        \begin{tikzpicture}
            \filldraw[fill=yellow!80!black,line width=0.5pt,fill opacity=0.6] (0,0) -- (3,0) -- (0,3) -- (0,0);
            \filldraw[fill=yellow!80!black,line width=0.5pt,fill opacity=0.6] (0,0) -- (3,0) -- (1.4,0.5) -- (0,0);
            \filldraw[fill=yellow!80!black,line width=0.5pt,fill opacity=0.6] (0,0) -- (0,3) -- (0.5,1.4) -- (0,0);
            \filldraw[fill=yellow!80!black,line width=0.5pt,fill opacity=0.6] (3,0) -- (0,3) -- (1.2,1.2) -- (3,0);
            \filldraw[fill=yellow!80!black,line width=0.5pt,fill opacity=0.6] (3,0) -- (1.4,0.5) -- (1.2,1.2) -- (3,0);
            \filldraw[fill=yellow!80!black,line width=0.5pt,fill opacity=0.6] (0,3) -- (0.5,1.4) -- (1.2,1.2) -- (0,3);
            \filldraw[fill=yellow!80!black,line width=0.5pt,fill opacity=0.6] (0,0) -- (0.5,1.4) -- (0.7,1) -- (0,0);
            \filldraw[fill=yellow!80!black,line width=0.5pt,fill opacity=0.6] (0,0) (0.7,1) -- (1.2,1.2) -- (0.5,1.4) -- (0,0);
            \coordinate[Bullet=black, label=below:1] (n1) at (0,0);
            \coordinate[Bullet=black, label=below:2] (n2) at (3,0);
            \coordinate[Bullet=black, label=left:3] (n3) at (0,3);
            \coordinate[Bullet=black, label=left:{4}] (n4) at (0.5,1.4);
            \coordinate[Bullet=black, label=below:{$n-1$}] (n5) at (1.4,0.5);
            \coordinate[Bullet=black, label=right:{$n$}] (n6) at (1.2,1.2);
            \coordinate[Bullet=black, label=below:{5}] (n7) at (0.7,1);
            \draw[dashed,red] (0.343,0) -- (1.042,1.958);
            \draw[dashed,red] (0,1.2) -- (2.52,0.48);
        \end{tikzpicture}
        \begin{tikzpicture}
            \filldraw[fill=yellow!80!black,line width=0.5pt,fill opacity=0.6] (0,0) -- (3,0) -- (0,3) -- (0,0);
            \filldraw[fill=yellow!80!black,line width=0.5pt,fill opacity=0.6] (0,0) -- (3,0) -- (1.4,0.5) -- (0,0);
            \filldraw[fill=yellow!80!black,line width=0.5pt,fill opacity=0.6] (0,0) -- (0,3) -- (0.5,1.4) -- (0,0);
            \filldraw[fill=yellow!80!black,line width=0.5pt,fill opacity=0.6] (3,0) -- (0,3) -- (1.2,1.2) -- (3,0);
            \filldraw[fill=yellow!80!black,line width=0.5pt,fill opacity=0.6] (3,0) -- (1.4,0.5) -- (1.2,1.2) -- (3,0);
            \filldraw[fill=yellow!80!black,line width=0.5pt,fill opacity=0.6] (0,3) -- (0.5,1.4) -- (1.2,1.2) -- (0,3);
            \filldraw[fill=yellow!80!black,line width=0.5pt,fill opacity=0.6] (0,0) -- (0.5,1.4) -- (0.7,1) -- (0,0);
            \filldraw[fill=yellow!80!black,line width=0.5pt,fill opacity=0.6] (0.5,1.4) -- (0.7,1) -- (1.2,1.2) -- (0.5,1.4);
            \filldraw[fill=yellow!80!black,line width=0.5pt,fill opacity=0.6] (0,0) -- (1.4,0.5) -- (1.2,1.2) -- (0.7,1) -- (0,0);
            \coordinate[Bullet=black, label=below:1] (n1) at (0,0);
            \coordinate[Bullet=black, label=below:2] (n2) at (3,0);
            \coordinate[Bullet=black, label=left:3] (n3) at (0,3);
            \coordinate[Bullet=black, label=left:{4}] (n4) at (0.5,1.4);
            \coordinate[Bullet=black, label=below:{$n-1$}] (n5) at (1.4,0.5);
            \coordinate[Bullet=black, label=right:{$n$}] (n6) at (1.2,1.2);
            \coordinate[Bullet=black, label=below:{5}] (n7) at (0.7,1);
            \draw[dashed] (n5) -- (n7);
        \end{tikzpicture}
        \caption{The first figure shows the lines to which $\overline{\mathbf{q}}(4),\overline{\mathbf{q}}(n-1),\overline{\mathbf{q}}(n)$ are constrained.\\
        The second figure shows how 2-simplices $12(n-1),34n$ introduce dependencies between the three variables $r,s,t$.\\
        The third figure shows how subsequently added equatorial vertices are constrained based on those previously added.\\
        The fourth figure shows a completed bipyramid.}
        \label{fig:bipcount}
\end{figure}

The following result demonstrates that, unlike in bar-joint rigidity \cite{connelly2005generic}, \cite{gortler2010characterizing}, global rigidity is not a generic property in volume rigidity.

\begin{corollary}\label{cor:notggr}
    The bipyramid $B_5$ admits rigid generic frameworks that are both globally rigid and not globally rigid.
\end{corollary}

\begin{proof}
    Recall from the proof of \cref{thm:bipcount} that the congruence classes of generic frameworks $(B_{n-2},\mathbf{p})$ correspond to solutions to the degree-$(n-4)$ univariate polynomial equation $f(\mathbf{p})(t)=0$.
    Moreover, $t=0$ is always a root, with the corresponding congruence class being $[\mathbf{p}]$.
    
    When $n=7$, $f(\mathbf{p})$ is cubic, and therefore has either one, two or three real solutions depending on whether the discriminant $\Delta(\mathbf{p}):=\Delta\left(\frac{1}{t}f(\mathbf{p})\right)$ is below, equal to or above zero respectively.
    Globally rigid frameworks are precisely those for which $\Delta(\mathbf{p})<0$.
    We note that, if $\mathbf{p}$ is generic, $\Delta(\mathbf{p})\neq0$, as the discriminant has coefficients in $\mathbb{Q}$.

    Below are the configuration matrices of the pinnings of two non-generic frameworks $(B_5,\mathbf{p}_1)$ and $(B_5,\mathbf{p}_2)$ with $\Delta(\mathbf{p}_1)<0$ and $\Delta(\mathbf{p}_2)>0$.
    Although they are not generic, since $\Delta$ is a continuous function, the positions of their unpinned vertices can be perturbed to yield pinnings of generic configurations $\mathbf{p}_1',\mathbf{p}_2'$ so that $\Delta(\mathbf{p}_1')<0$ and $\Delta(\mathbf{p}_2')>0$.
    \begin{equation*}
    \begin{split}
        C(\overline{\mathbf{p}_1}) &= \begin{bmatrix}
            1 & 1 & 1 & 1 & 1 & 1 & 1 \\
            0 & 1 & 0 & \frac{1}{5} & \frac{1}{7} & \frac{1}{11} & \frac{1}{2} \\
            0 & 0 & 1 & \frac{1}{13} & \frac{1}{19} & \frac{1}{17} & \frac{1}{2}
        \end{bmatrix}, \\
        C(\overline{\mathbf{p}_2}) &= \begin{bmatrix}
            1 & 1 & 1 & 1 & 1 & 1 & 1 \\
            0 & 1 & 0 & \frac{1}{7} & \frac{1}{5} & \frac{1}{41} & \frac{1}{2} \\
            0 & 0 & 1 & \frac{1}{19} & \frac{1}{17} & \frac{1}{13} & 20
        \end{bmatrix}.
    \end{split}
    \end{equation*}
\end{proof}

\bibliographystyle{plainurl}
\bibliography{BoundsTrisS2}

\end{document}